\newfont{\bcb}{msbm10}
\newfont{\matb}{cmbx10}
\newfont{\got}{eufm10}
\documentclass{article}
\usepackage[cp1250]{inputenc}
\usepackage{amsmath}
\usepackage{amsmath,amsthm}
\usepackage{amssymb,latexsym}
\usepackage{latexsym}
      \usepackage[all, knot]{xy}
     \xyoption{arc}
\usepackage{hyperref}

\DeclareMathOperator{\tr}{tr}
\DeclareMathOperator{\rank}{rank}

\DeclareMathOperator{\extr}{extr}
\providecommand{\abs}[1]{\lvert#1\rvert}
\providecommand{\norm}[1]{\lVert#1\rVert}

\begin{document}
\title{Note on the spread of real symmetric matrices with entries in fixed interval.}
\author{Iwo Biborski(*)}

\footnotetext{(*) Authors e-mail address: iwo.biborski@gmail.com}

\maketitle

\begin{abstract}
The spread of a matrix is defined as the maximum of distances between any two eigenvalues of that matrix. In this paper we investigate spread maximization as a function on compact convex subset of the set of real symmetric matrices. We provide some general results and further, we study spread maximizing problem on $S_{n}\left[a,b\right]$(the set of symmetric matrices with entries restricted to the interval $\left[a,b\right]$). In particular, we develop some results by X. Zhan (see \cite{[Zhan]}), S. M. Fallat and J. J. Xing (see \cite{[FaXi]}).

\end{abstract}

\section{Introduction}
Let $A$ be $n\times n$ matrix over field $\mathbb{C}$ or $\mathbb{R}$, and let $\lambda_{1},\dots,\lambda_{n}$ be its eigenvalues. Following by \cite{[Mir]}, we recall a definition of the spread:

\begin{equation}\tag{1}
s_{n}\left(A\right) =\max\{\abs{\lambda_{i} - \lambda_{j}}\}
\end{equation}
for $i,j\in{1,\dots,n}$. If the spectrum of $A$ consists only real numbers, we can assume the following order of the eigenvalues
$$\lambda_{1}\geq\lambda_{2}\geq\dots\geq\lambda_{n},$$ 
and in that case

\begin{equation}\tag{2}
s_{n}\left(A\right) =\lambda_{1} - \lambda_{n}.
\end{equation}

Straightforward from the definition of the spread, for all $A,B\in S_{n}(\mathbb{R})$, we have
\begin{gather*}s_{n}(\alpha A) = \abs{\alpha}s_{n}(A), \tag{3}\\ 
\forall_{\alpha \in \,\left[0,1\right]} : s_{n}(\alpha A + (1-\alpha) B)\leq \alpha s_{n}(A) + (1-\alpha)s_{n}(B),\tag{4}
\end{gather*}
in particular $s_{n}$ is a convex function (for ref. see for instance \cite{[Merik2]}). Through the paper, by $S_{n}(\mathbb{R})$ we denote the set of $n\times n$ real symmetric matrices. Further, by $S_{n}\left[a,b\right]$ we mean the subset of $S_{n}(\mathbb{R})$ of matrices with entries from an interval $\left[a,b\right]$. We always assume that eigenvector is unitary.

The spread of a matrix is a part of width study of characteristic polynomial roots locus. It also plays role in the combinatorial optimization theory (see \cite{[FiBuRe]}) and graph theory (see \cite{[GreHerKi]} for reference). One of the challenging problems is to give possibly best bounds of the spread, both in general and for some special types of matrices as well. Through the years a lot of upper and lower bounds of the spread where given, for reference for the upper bounds see \cite{[Mir]}, \cite{[Mir2]}, \cite{[De]}, \cite{[Bes]}, and more recent papers such \cite{[Wu]}, \cite{[ShaKu]} or \cite{[Ping]}. For the lower bounds see \cite{[Jon]} and \cite{[Merik]}.

In this paper we investigate the maximal value of the spread restricted to a convex compact subset $V$ of the set of real symmetric matrices. Let us recall that in the papers \cite{[Zhan]} and \cite{[FaXi]} the authors, among the other results, show that the spread restricted to $V = S_{n}\left[a,b\right]$ is maximized by some matrix which entries belong to $\{a,b\}$. Our research is motivated by the following problem: does exist a matrix from $S_{n}\left[a,b\right]$ which maximizes the spread on this set and is not extreme point of $S_{n}\left[a,b\right]$? We give a negative answer to that question in \textbf{Theorem 5}. However, to achieve that goal, we investigate more general problem: under which hypothesis on $V$, the spread restricted to $V$ is maximized only by some extreme points? We give necessary and sufficient condition for $V$ to have that property and use this characterization to prove \textbf{Theorem 5}. In the last section we provide all maximization matrices of $s_{3}$ restricted to $S_{3}\left[0,1\right]$.

\section{The maximal spread on the convex compact subsets of the set of real symmetric matrices.}
In this section we consider the following problem: let $V$ be a convex compact subset of $S_{n}(\mathbb{R})$. Under what conditions, $\max\{s_{n}(A): A\in V\}$ is attained only by extreme points of $V$? 

Let us recall that for a convex set $C$, an element $c\in C$ is extreme point if and only if $c$ does not belong to the interior of any non-degenerated interval contained in $C$. In other words, $c$ is an extreme point of $C$ if, for all $a,b \in C$ and $t\in\left[0,1\right]$ the equality $ta+(1-t)b = c$ implies that $t=0$ or $t=1$. By $\extr(C)$ we denote the set of all extreme points of $C$. 

Let $V$  be a convex compact subset of $S_{n}(\mathbb{R})$. We define
\begin{equation}
V_{max}^{s_{n}} := \{A\in V: s_{n}(A) = \max\{s_{n}(B):\,B\in V\}\}. \tag{5}
\end{equation}
For $A\in S_{n}(\mathbb{R})$, by $\lambda_{1}^{A}$ we denote maximal, and by $\lambda_{n}^{A}$ minimal eigenvalue of $A$. Further, we define 
\begin{gather*}
u(A) := \{u\in \mathbb{R}^{n}: Au = \lambda_{1}^{A}u \},\\
v(A) := \{v\in \mathbb{R}^{n}: Av = \lambda_{n}^{A}v \},
\end{gather*}
for the set of eigenvectors corresponding with maximal and minimal eigenvalue respectively.

\vspace{2ex}
\textbf{Definition 1. Extreme spread maximization.}
\textit{We say that $V$ has the property of extreme spread maximization if $V_{max}^{s_{n}}\subset \extr(V)$.}

\vspace{2ex}
We give two examples of sets which does not have the property of extreme spread maximization.

\vspace{2ex}
\textit{Example 1.} The set of diagonal matrices which are multiplication of identity matrix $I$ by real number from closed interval doses not have the property of extreme spread maximization. It follows from the fact that for all matrices from that set, the spread is constantly equal $0$. 

\vspace{2ex}
\textit{Example 2.} Let $W$ be a convex compact subset of $S_{n}(\mathbb{R})$ which has extreme spread maximization property. Consider the following set:
$$W_{a,b} := \{X:\, X = Y+\theta I,\,Y\in W,\,\,\theta\in\left[a,b\right]\}.$$
If $X$ maximizes the spread on $W$, thus for each $\theta\in\left[a,b\right]$, $X+\theta I$ maximizes the spread on $W_{a,b}$, thus $W_{a,b}$ does not have the property from the \textbf{Definition 1}.

\vspace{2ex}
The goal of this section is to give necessary and sufficient condition for $V$ to have property of extreme spread maximization. To achieve this, we need two following lemmas:

\vspace{2ex}
\textbf{Lemma 1.} \textit{Let} $A,B\in S_{n}(\mathbb{R})$. Then $s_{n}(A+B)=s_{n}(A)+s_{n}(B)$ \textit{if and only if} $u(A)\cap u(B)\neq\emptyset$ and $v(A)\cap v(B)\neq\emptyset$

\begin{proof}
Suppose that  $u(A)\cap u(B)\neq\emptyset$ and $v(A)\cap v(B)\neq\emptyset$. Let $u$ be the eigenvector corresponding with $\lambda_{1}^{A}$ and $\lambda_{1}^{B}$, and $v$ be the eigenvector corresponding with $\lambda_{n}^{A}$ and $\lambda_{n}^{B}$. Then we have
\begin{gather*}
s_{n}(A) + s_{n}(B) = \lambda_{1}^{A} - \lambda_{n}^{A} + \lambda_{1}^{B} - \lambda_{n}^{B} = \\
u^{T}Au - v^{T}Av + u^{T}Bu - v^{T}Bv = u^{T}(A+B)u - v^{T}(A+B)v\leq s_{n}(A+B).
\end{gather*}
On the other hand, $s_{n}(A+B) \leq s_{n}(A) + s_{n}(B)$ (see \cite{[Merik2]}, COROLLARY 2 for more general result) and thus we obtain the equality.

\vspace{2ex}
Now, let us assume that $s_{n}(A+B)=s_{n}(A)+s_{n}(B)$. By THEOREM 1 from \cite{[Merik2]}, we have
\begin{equation}\tag{6}
\lambda_{1}^{A+B}\leq\lambda_{1}^{A}+\lambda_{1}^{B}
\end{equation}
and 
\begin{equation}\tag{7}
-\lambda_{n}^{A+B}\leq\-\lambda_{n}^{A}-\lambda_{n}^{B}.
\end{equation}
Therefore, according to (6) and (7),
\begin{equation*}
s_{n}(A+B)=s_{n}(A)+s_{n}(B)
\end{equation*}
if and only if 
\begin{gather*}
\lambda_{1}^{A+B}=\lambda_{1}^{A}+\lambda_{1}^{B}
\end{gather*}
and
\begin{gather*}
    \lambda_{n}^{A+B}=\lambda_{n}^{A}+\lambda_{n}^{B}.
\end{gather*}
We will show that $u_{A+B}\in u(A+B)$ belongs to $u(A)\cap u(B)$.
For the contrary, suppose that for all $u_{A+B}\in u(A+B)$, $u_{A+B}\notin u(A)$. Then 
$$
    \lambda_{1}^{A+B} = u_{A+B}^{T}(A+B)u_{A+B}=u_{A+B}^{T}Au_{A+B}+u_{A+B}^{T}Bu_{A+B}<\lambda_{1}^{A}+\lambda_{1}^{B},
$$
and therefore $s_{n}(A+B)<s_{n}(A)+s_{n}(B)$-- a contradiction. In the similar way we prove that $u_{A+B}$ is an eigenvector of $B$ corresponding with $\lambda^{B}_{1}$. 

Analogical reasoning for $v_{A+B}$, $v_{A}$ and $v_{b}$ proves that $v_{A+B}$ is a common eigenvector for minimal eigenvalue of $A$ and $B$.
\end{proof}

\textbf{Lemma 2.} \textit{Let $f:X\longrightarrow \mathbb{R}$ be a convex function on a compact convex set $X$ such that, for all $t\in\left[0,1\right]$, $f(tx) = tf(x)$. Suppose that $f$ attains it's maximum at non-extreme point $c\in X$. Then, for all $a,\,b\in X$, if $ta+(1-t)b=c$ for some $t\in (0,1)$ we have }
\begin{equation}
f(a)=f(b)=f(c). \tag{8}
\end{equation}

\begin{proof}
Suppose that $f(a)<f(c)$ or $f(b)<f(c)$. Then, since $f$ is convex, we get
$$f(c)\leq t f(a)+(1-t)f(b)<t f(c)+(1-t)f(c)=f(c)$$
which is a contradiction. Therefore $f(a)\geq f(c)$ and $f(b)\geq f(c)$. Since $f(c)$ is maximal value, we obtain $f(a)=f(b)=f(c).$
\end{proof}

By (3) and (4), we see that $s_{n}$ fulfills the assumptions of \textbf{Lemma 2}. This leads us to the straightforward corollary.

\vspace{2ex}
\textbf{Corollary 3.}\textit{ Suppose that $V$ does not have the property of extreme spread maximization. Let $\alpha A + (1-\alpha)B\in V_{max}^{s_{n}}$ for some $A\neq B$ from $V$. Then,
\begin{gather*}
u(A)\cap u(B)\neq\emptyset, \\ v(A)\cap v(B)\neq\emptyset.
\end{gather*}
In particular $A$ and $B$ maximizes $s_{n}$ on $V$.}

\begin{proof}
Directly from \textbf{Lemma 2} we get $s_{n}(A)=s_{n}(B)=s_{n}(C)$. On the other hand, $$s_{n}(C)=s_{n}(\alpha A + (1-\alpha)B) = s_{n}(\alpha A) + s_{n}((1-\alpha)B).$$
Let us notice that for any matrix $X\in S_{n}(\mathbb{R})$ and real number $\gamma\geq 0$, 
\begin{gather*}
u(\gamma X) = u(X)\,\,\, \text{and} \,\,\,v(\gamma X)=v(X).
\end{gather*}
Therefore, by \textbf{Lemma 1}, we conclude that $$u(A)\cap u(B)\neq\emptyset$$ and $$v(A)\cap v(B)\neq\emptyset.$$
\end{proof}

\section{The maximal spread on \texorpdfstring{$S_{n}\left[a,b\right]$}{S{n}[a,b]}.}

Let $a,b\in\mathbb{R}$ and $a < b$. Our goal is to prove that if $A$ maximizes $s_{n}$ on $S_{n}\left[a,b\right]$, then $A$ is an extreme point of $S_{n}\left[a,b\right]$. Let us point out that extreme points of $S_{n}\left[a,b\right]$ are those matrices which entries belong to $\{a,b\}$. It is known fact that there exists extreme point from $S_{n}\left[a,b\right]$ which maximizes $s_{n}$ on that set, see \cite{[Zhan]}(proof of Theorem 1 and Theorem 7) and \cite{[FaXi]}(LEMMA 2.2) for reference.

\vspace{2ex}
Let us start with the following lemma:

\vspace{2ex}
\textbf{Lemma 4.}\textit{ For all $n\in\mathbb{N}_{+}$} 
$$\max\{s_{n}(A):\,A\in S_{n}\left[a,b\right]\}\,<\,\max\{s_{n+1}(A):\,A\in S_{n+1}\left[a,b\right]\}.$$

\begin{proof}
Let $A$ be a matrix for which $s_{n}$ attains maximum. Let $u$ and $v$ be the eigenvectors of $A$ corresponding with $\lambda_{1}^{A}$ and $\lambda_{n}^{A}$ respectively. We define a block matrix:
$$\hat{A}(X):=\left[\begin{array}{cc}a&X^{T}\\X&A\end{array}\right],$$
where $X\in [a,b]^{n}$. Put

\begin{gather*}\hat{u}:=\left[\begin{array}{c}0\\u\end{array}\right], \, \hat{v} :=\left[\begin{array}{c}0\\v\end{array}\right]
\end{gather*}
It is clear that 
\begin{equation*}\hat{u}^{T}\hat{A}(X)\hat{u} - \hat{v}^{T}\hat{A}(X)\hat{v} = u^{T}Au- v^{T}Av = s_{n}(A),
\end{equation*}
thus $s_{n}(A) \leq s_{n+1}(\hat{A}(X))$. In particular 
\begin{align*}\tag{9}
\lambda_{1}^{\hat{A}(X)} \geq \lambda_{1}^{A}, \\ 
\lambda_{n}^{\hat{A}(X)} \leq \lambda_{n}^{A}.
\end{align*}

Suppose, that for all $X\in [a,b]^{n}$, $s_{n+1}(\hat{A}(X)) = s_{n}(A)$.
Then $\hat{A}(X)$ has no greater eigenvalue than $\lambda_{1}^{A}$ and no smaller than $\lambda_{n}^{A}$ and, thanks to (9), we have
\begin{gather*}
\lambda_{1}^{\hat{A}(X)} = \lambda_{1}^{A}, \\
\lambda_{n}^{\hat{A}(X)} = \lambda_{n}^{A}, \\
\end{gather*}
for all $X$. In addition, $\hat{u}$ must be the eigenvector of $\hat{A}(X)$ corresponding with $\lambda_{1}^{\hat{A}(X)}$. Since $u$ is non-zero vector, we can choose $X_{0}\in [a,b]^{n}$ such that $X_{0}^{T}u\neq 0$. Therefore
$$\hat{A}(X_{0})\hat{u} =\left[\begin{array}{cc}a & X_{0}^{T}\\X_{0} & A\end{array}\right]\left[\begin{array}{c}0\\u\end{array}\right] = \left[\begin{array}{c}X_{0}^{T}u\\ \lambda_{1}^{A}u\end{array}\right].$$ On the other hand $\hat{A}(X_{0})\hat{u}=\left[\begin{array}{c}0\\ \lambda_{1}^{A}u\end{array}\right]$, which contradicts with $X_{0}^{T}u\neq 0$.
\end{proof}

Now we are ready to formulate and prove our main result which we present below:

\vspace{2ex}
\textbf{Theorem 5.}\textit{ Let $C$ be a matrix which maximizes $s_{n}:S_{n}\left[a,b\right]\rightarrow\mathbb{R}$. Then $C\in\extr(S_{n}\left[a,b\right])$.}

\begin{proof}
For the contrary, suppose that $C$ maximize $s_{n}$ but is non-extreme point of $S_{n}\left[a,b\right]$. Let $c_{ij}$ be the entry of $C$ which belongs to $(a,b)$. Let $A = (a_{kl})_{k,l},\, B=(b_{kl}){k,l}\in S_{n}\left[a,b\right]$ be the symmetric matrices which entries are defined as follows:
\begin{equation}\tag{10}
a_{kl}=\begin{cases}
c_{kl}, & k\neq i \lor l\neq j\\
a, & k=i \land l=j,
\end{cases}
\end{equation}
\begin{equation}\tag{11}
b_{kl}=\begin{cases}
        c_{kl}, & k\neq i \lor l\neq j\\
        b, & k=i \land l=j.
        \end{cases}
\end{equation}
Since $c_{ij} = \alpha a + (1-\alpha)b$ for some $\alpha\in\left(0,1\right )$, thus $C$ is a convex combination of $A$ and $B$ and $C = \alpha A + (1 - \alpha)B$. By \textbf{Lemma 2},
\begin{gather*}s_{n}(C) = s_{n}(A) = s_{n}(B).
\end{gather*}
Further, from \textbf{Corollary 3}, there exist $u$ and $v$ such that 
\begin{gather*}
\begin{cases}
    Au=\lambda_{1}^{A} u \\ Bu=\lambda_{1}^{B}u \\Av=\lambda_{n}^{A}v \\Bv=\lambda_{n}^{B}v. 
\end{cases}
\end{gather*}
Let $X = B-A$. Then, from (10) and (11), $X$ has all entries equal zero except coefficients $(i,j)$ and $(j,i)$ which are equal $b-a$. Since

$$\lambda^{B}_{1}u = (A +X)u = Au + Xu = \lambda^{A}_{1}u +Xu,$$ we obtain
$$(\lambda^{B}_{1}-\lambda^{A}_{1})u = Xu.$$
We can write previous equality in the following vector form:
\begin{equation}\tag{12}
\left[\begin{array}{c}(\lambda^{B}_{1}-\lambda^{A}_{1})u_{1} \\ \vdots \\ (\lambda^{B}_{1}-\lambda^{A}_{1})u_{i} \\ \vdots \\ (\lambda^{B}_{1} - \lambda^{A}_{1})u_{j} \\ \vdots \\ (\lambda^{B}_{1}-\lambda^{A}_{1})u_{n}\end{array}\right]=\left[\begin{array}{c}0 \\ \vdots \\ (b-a)u_{j} \\ \vdots \\ (b-a)u_{i} \\ \vdots \\ 0 \end{array}\right].
\end{equation}
In the same fashion, we get equality $v$ and the eigenvalues $\lambda_{n}^{A}$ and $\lambda_{n}^{B}$:
\begin{equation}\tag{13}\left[\begin{array}{c}(\lambda^{B}_{n}-\lambda^{A}_{n})v_{1} \\ \vdots \\ (\lambda^{B}_{n}-\lambda^{A}_{n})v_{i} \\ \vdots \\ (\lambda^{B}_{n}-\lambda^{A}_{n})v_{j} \\ \vdots \\ (\lambda^{B}_{n} -\lambda^{A}_{n})v_{n}\end{array}\right]=\left[\begin{array}{c}0  \vdots \\ (b-a)v_{j} \\ \vdots \\ (b-a)v_{i} \\ \vdots \\ 0\end{array}\right].\end{equation}
Since $s_{n}(A) = s_{n}(B)$, it is clear that $\lambda_{1}^{A} = \lambda_{1}^{B}$ if and only if $\lambda_{n}^{A} = \lambda_{n}^{B}$. Comparing corresponding entries of vectors from (12) and (13), we conclude that either 
\begin{equation}\tag{14}
\begin{cases}
    \lambda_{1}^{A} = \lambda_{1}^{B} \\
    \lambda_{n}^{A} = \lambda_{n}^{B}
\end{cases}
\end{equation}
or 
\begin{equation}\tag{15}
u_{k} = v_{k} = 0\,\,\, \text{for}\,\,\, k\neq i \,\,\, \text{and}\,\,\, k\neq j.
\end{equation}

If (14) holds, then
\begin{gather*}
    \begin{cases}
      u_{i}=u_{j}= 0\\
    v_{i}=v_{j}=0
\end{cases}
\end{gather*}
and $s_{n}(A)$ is equal to the spread of a matrix obtained by crossing out $i$-th and $j$-th rows and columns from $A$. Therefore
\begin{gather*}
\max\{s_{n}(A):\,A\in S_{n}\left[a,b\right]\}\,\leq\,\max\{s_{n-2}(A):\,A\in S_{n-2}\left[a,b\right]\},
\end{gather*}
or
\begin{gather*}
\max\{s_{n}(A):\,A\in S_{n}\left[a,b\right]\}\,\leq\,\max\{s_{n-1}(A):\,A\in S_{n-2}\left[a,b\right]\},
\end{gather*}
if $i=j$.
This contradicts with \textbf{Lemma 4}. 

If (15) is true, $s_{n}(A)$ is equal to the spread of $2\times2$ matrix obtained by crossing out all rows and columns except $i$-th and $j$-th, or $1\times 1$ matrix if $i=j$. In that case  
\begin{gather*}
\max\{s_{n}(A):\,A\in S_{n}\left[a,b\right]\}\,\leq\,\max\{s_{2}(A):\,A\in S_{2}\left[a,b\right]\},
\end{gather*}
or
\begin{gather*}
\max\{s_{n}(A):\,A\in S_{n}\left[a,b\right]\}\,\leq\,\max\{s_{1}(A):\,A\in S_{2}\left[a,b\right]\}.
\end{gather*}
This contradicts with \textbf{Lemma 4} as well.
\end{proof}

Let us recall that $s_{n}(A)$ can be expressed in the following form:
\begin{equation}\tag{16}
s_{n}(A) = e \left[A\circ(uu^{T} - vv^{T})\right]e^{T}
\end{equation}
where $e = \left[1,1,\dots,1\right]$, $u$ and $v$ are eigenvectors corresponding with greatest and smallest eigenvalue respectively, and "$\circ$" stands for the Hadamard product (see \cite{[Zhan]} and \cite{[FaXi]}). This representation was used in \cite{[FaXi]} to prove the following fact: for $a\in\left[-1,1\right)$, if $A\in S_{n}\left[a,1\right]$ is of rank 2 and maximizes $s_{n}$, then it must be in some special block form and 
\begin{equation}\tag{17}
u_{i}u_{j} - v_{i}v_{j}\neq 0
\end{equation} for each $(i,j)$, (\cite{[FaXi]}, THEOREM 2.5 and 2.6). Here we can use our \textbf{Theorem 5} to go further and prove that (17) holds for $A$ which maximizes $s_{n}$ on $S_{n}\left[a,b\right]$ without any additional assumptions on form and rank. We have the following

\vspace{2ex}
\textbf{Corollary 6.}\textit{ Let $A=(a_{ij})_{i,j}$ be a matrix which maximizes $s_{n}$ on $S_{n}\left[a,b\right]$ and let $u$ and $v$ be eigenvectors corresponding with maximal and minimal eigenvalue, respectively. Then $u_{i}u_{j} - v_{i}v_{j}\neq 0$ for each $i,j\in \{1, \dots n\}$.}
\begin{proof}
Suppose that there exist $i$ and $j$ such that $u_{i}u_{j} - v_{i}v_{j} = 0$. Then, by (16), $s_{n}(A)$ does not depend on $a_{ij}$. Put $\tilde{A}$ for matrix which all entries are the same as for $A$ except $(i,j)$-th and $(j,i)$-th elements, where we place some element from $(a,b)$. Then $s_{n}(\tilde{A}) = s_{n}(A)$ which contradicts with \textbf{Theorem 5}, since $\tilde{A}$ is non-extreme point of $S_{n}\left[a,b\right]$. 
\end{proof}

Another consequence of \textbf{Theorem 5} is a straightforward condition on the diagonal of maximizing matrix:

\vspace{2ex}
\textbf{Corollary 7.} \textit{Let $A\in S_{n}\left[a,b\right]$ be the matrix that maximizes $s_{n}$ on  $S_{n}\left[a,b\right]$. Then the diagonal of $A$ must contain at least one element equal $a$ and at least one element equal $b$.}

\vspace{2ex}
\begin{proof}
It follows from \textbf{Theorem 5} that the diagonal of $A$ has entries belonging to $\{a,b\}$. Suppose that all entries of the diagonal are equal. If all entries are equal $a$, we take $\gamma > 0$ that $a+ \gamma < b$. If each entry of the diagonal is $b$, we take $\gamma < 0$ such that $b + \gamma > a$. Then, since $s_{n}(A + \gamma I) = s_{n}(A)$, we conclude that $A+\gamma I$ maximizes $s_{n}$ on $S_{n}\left[a,b\right]$. This contradicts with \textbf{Theorem 5}.  
\end{proof}

\section{The maximal spread on \texorpdfstring{$S_{3}\left[0,1\right]$}{S{3}[0,1]}.}

In this section we will prove that matrices
$$Y := \left[\begin{array}{ccc}0&1&1\\ 1&1&1\\ 1&1&1 \end{array}\right],\,Y_{1} := \left[\begin{array}{ccc}1&1&1\\ 1&0&1\\ 1&1&1 \end{array}\right],\,Y_{2} := \left[\begin{array}{ccc}1&1&1\\ 1&1&1\\ 1&1&0 \end{array}\right]$$
maximize the spread on $S_{3}\left[0,1\right]$ and that there are no other maximizing matrices on that set. Since $Y_{1}$ and $Y_{2}$ are similar to $Y$, it is enough to prove that $s_{3}(Y)$ is maximal on $S_{3}\left[0,1\right]$. By simple computation $s_{3}(Y) = 2\sqrt{3}$. Before we present the proof, let us recall well known upper bound of the spread due to Mirsky (see \cite{[Mir]}):
\begin{equation}\tag{18}
s_{n}(A) \leq \sqrt{2\norm{A}_{F}^{2}-\frac{\tr(A)^{2}}{n}},
\end{equation}
where $\norm{A}_{F}$ is the Frobenius norm. We will use an inequality (18) to prove the following 

\vspace{2ex}
\textbf{Corollary 8.} \textit{Matrix $Y$ maximizes the spread on $S_{3}\left[0,1\right]$ and, up to permutation of eigen basis, is the only maximizing matrix on that set.}

\begin{proof}
Let $X$ be a matrix that maximizes $s_{3}$ on $S_{3}\left[0,1\right]$. By \textbf{Theorem 5}, $X$ is extreme point of $S_{3}\left[0,1\right]$ so all entries of $X$ belong to $\{0,1\}$. By THEOREM 2.5 from \cite{[FaXi]}, $Y$ maximizes the spread for all rank two matrices from $S_{3}\left[0,1\right]$. Since the greatest spread of matrix with rank one is equal $3$, is enough to show that there is no matrix of rank $3$ from $S_{3}\left[0,1\right]$ which maximizes the spread. 

For the contrary, suppose that there is $X\in S_{3}\left[0,1\right]$ with rank $3$, that maximizes $s_{3}$. By \textbf{Theorem 5} $X$ is an extreme point of $S_{3}\left[0,1\right]$, thus it must have at least two zeros among its coefficients and rest of them equal $1$. We will consider two cases - first, where $X$ has two zeros and then $3$ or more zeros. 

\textit{Case 1.} We assume that $X$ has two zeros, thus is equal to

$$X_{1} = \left[\begin{array}{ccc}0&1&1\\ 1&0&1\\ 1&1&1 \end{array}\right],$$

or

$$X_{2} = \left[\begin{array}{ccc}1&0&1\\ 0&1&1\\ 1&1&1 \end{array}\right],$$
or any other matrix obtained from $X_{1}$ or $X_{2}$ by changing an ordering of basis. By \textbf{Corollary 7}, $X_{2}$ cannot maximize $s_{3}$. On the other hand, 
$$s_{3}(X_{1}) = 2+\sqrt{2} < 2\sqrt{3} = s_{3}(Y),$$
therefore $X_{1}$ does not maximize the spread as well.

\textit{Case 2.} Now we assume that $X$ is an extreme point of $S_{3}\left[0,1\right]$ with $3$ or more zeros. It implies that $\norm{X}_{F}^{2} \leq 6$. But according to (18) we have

$$s_{3}(X) \leq \sqrt{12-\frac{\tr(X)^{2}}{3}} \leq 2\sqrt{3}.$$
By \textbf{Corollary 7}, if $X$ maximizes $s_{3}$ on $S_{3}\left[0,1\right]$ then $\tr(A)$ must be greater than zero. Therefore
\begin{gather*}
    s_{3}(X) < 2\sqrt{3}.
\end{gather*}
This completes our proof.
\end{proof}
Thanks to \textbf{Corollary 8} we give a positive answer to the conjecture by Fallat and Xing for special case $n=3$. We recall the conjecture below:

\vspace{2ex}
\textbf{Conjecture (\cite{[FaXi]}, CONJECTURE 2.7)}\textit{ For $A\in S_{n}\left[a,b\right]$, the spread in this class is attained by $A\in  S_{n}\{a,b\}$ and $\rank(A) = 2$.}

\vspace{2ex}
Last conjecture remains open in general. However our results give more precise information about spread maximization on $S_{n}\left[a,b\right]$ which enhance to the further investigation.

\end{document}